\documentclass[11pt,a4]{article}

\usepackage[T2A]{fontenc}
\usepackage[utf8]{inputenc}
\usepackage{amsfonts, amsmath, amsthm}
\usepackage{amssymb}

\newcommand{\R}{{\mathbb R}}
\newcommand{\E}{{\mathbb E}}
\newcommand{\Prob}{{\mathbb P}}
\def\e{\varepsilon}

\def\myproof{{\it Proof}  }

\newtheorem{Prop}{Proposition}[section]

\newtheorem{Cor}{Corollary}[section]

\newtheorem{Th}{Theorem}[section]

\begin{document}
%
%
%


\begin{center}
\textbf{Parallel algorithms and probability of large deviation for stochastic optimization problems}

\end{center}

\begin{center}
\textit{Pavel Dvurechensky
\footnote{corresponding author, pavel.dvurechensky@wias-berlin.de, Weierstrass Institute for Applied Analysis and Stochastics, Mohrenstr. 39, 10117 Berlin, Germany; Institute for Information Transmission Problems RAS, Bolshoy Karetny per. 19, build.1,
Moscow, Russia 127051;\\
}, 
Alexander Gasnikov
\footnote{gasnikov.av@mipt.ru, Moscow Institute of Physics and Technology, 9 Institutskiy per., Dolgoprudny, Moscow Region,
Russia 141700; Institute for Information Transmission Problems RAS, Bolshoy Karetny per. 19, build.1,
Moscow, Russia 127051; \\ 
}, 
Anastasia Lagunovskaya
\footnote{a.lagunovskaya@phystech.edu, Moscow Institute of Physics and Technology, 9 Institutskiy per., Dolgoprudny, Moscow Region,
Russia 141700;\\
} 
}
\end{center}

\begin{abstract}
We consider convex stochastic optimization problems under different assumptions on the properties of available stochastic subgradient. It is known that, if the value of the objective function is available, one can obtain, in parallel, several independent approximate solutions in terms of the objective residual expectation. 
Then, choosing the solution with the minimum function value, one can control the probability of large deviation of the objective residual. On the contrary, in this short paper, we address the situation, when the value of the objective function is unavailable or is too expensive to calculate. 
Under "`light-tail"' assumption for stochastic subgradient and in general case with moderate large deviation probability,  we show that parallelization combined with averaging gives bounds for probability of large deviation similar to a serial method. Thus, in these cases, one can benefit from parallel computations and reduce the computational time without loss in the solution quality.

\end{abstract}





\section{Introduction}
We consider the following general stochastic optimization problem over a convex compact set $Q$:
\begin{equation}
\label{eq1}
\min_{x\in Q \subset E} \left\{f\left( x \right):=\E_\xi \left[ {f\left( {x,\xi } \right)} \right]\right\},
\end{equation}
where $E$ is a finite-dimensional real vector space, $\xi$ is a random vector, $f\left( {x,\xi } \right)$ is a closed convex function w.r.t $x$ for a.e. $\xi$, and $x$ and $\xi$ are independent. Under these assumptions, this problem is a convex optimization problem.

Our main goal is to approximately solve this problem using some algorithm. Usually, in the stochastic optimization literature \cite{1, 2}, two measures for the quality of an approximate solution $\bar{x}$ are considered. The first is expectation of the objective residual. In this case, $\bar{x}$ is an $\e$-solution of \eqref{eq1} for $\e >0$ iff  $\E f(\bar{x}) - f_* \leq \e$, where $f_*$ is the optimal value in \eqref{eq1} and expectation is taken with respect to all the randomness arising in the algorithmic process. The second is bound for probability of large deviation of the objective residual. In this case $\bar{x}$, is an $(\e,\sigma)$-solution of \eqref{eq1} for $\e >0$, $\sigma \in (0,1)$ iff  $\Prob\{f(\bar{x}) - f_* > \e \} \leq \sigma$. We mainly focus on the latter quality measure in this paper. 

It is known that, if the value of the objective function is available (e.g. in randomized methods \cite{3}), one can obtain in parallel logarithmic in $\sigma^{-1}$ number of independent $\e$-solutions. Then, the solution with the minimum function value is an $(\e,\sigma)$-solution. Nevertheless, the value of the objective function can be unavailable or too expensive to calculate or approximate. The latter is easy to imagine since, for $\xi \in \R^p$, the computational effort for calculation of $\bar{f}$ s.t. $| \bar{f} - f(x)| \leq \delta$ can amount up to $O\left(\delta^{-p}\right)$ calculations of $f(x,\xi)$ at different $\xi$. 
Our goal is to propose a technique, which allows to obtain an $(\e,\sigma)$-solution based on a number computed in parallel $\e$-solutions without calculation of the function $f(x)$ value.

Our approach is based on Stochastic Mirror Descent algorithm \cite{2}. It turns out that, under some mild assumptions, an $\e$-solution to \eqref{eq1}, obtained by Stochastic Mirror Descent, is also an $(\tilde{\e},\sigma)$-solution. We use this fact and calculate in parallel logarithmic in $\sigma^{-1}$ number of independent $\e$-solutions, average them and prove that this average is an $(\e,\sigma)$-solution. Thus, we can benefit from parallelization and reduce the computation time without any loss in the solution quality.


\section{Stochastic Mirror Descent}
This section is devoted to description of Stochastic Mirror Descent (SMD) \cite{2,3} and its convergence properties in terms of expectation of the objective residual and also in terms of probability of large deviation of this residual. These convergence results provide the basis of our approach for constructing an $(\e,\sigma)$-solution by parallelization.
Let us choose some norm $\left\|\cdot \right\|$ on $E$ and denote the conjugate norm by $\left\| \cdot \right\|_\ast $.
We assume that, at any point $x\in Q$, a stochastic subgradient $\nabla _x f\left( {x,\xi } \right)$ of $f(x)$ is available and satisfies
\begin{equation}
\E_\xi \left[ {\nabla _x f\left( {x,\xi } \right)} \right] \in \partial f\left( 
x \right),
\quad
\E_\xi \left[ {\left\| {\nabla _x f\left( {x,\xi } \right)} \right\|_\ast 
^2 } \right]\le M^2
\label{eq:As}
\end{equation}
for some constant $M >0$.
We choose a prox-function $d\left( x \right)$ which is 
1-strongly convex in $\left\| \cdot \right\|$.
Let $x^0 = \arg \min_{x\in Q} d(x)$. W.l.o.g. we assume that $d(x^0)=0$. The algorithm uses Bregman's divergence
$V_z\left( {x} \right)=d\left( x \right)-d\left( z \right)-\left\langle {\nabla d\left( z \right),x-z} \right\rangle$.
Let $x_*$ be a solution of \eqref{eq1}, $R$ be a number s.t. $V_{x^0}\left( {x_\ast} \right) \leq R^2$, and $\bar{R}$ be a number s.t. $\max_{x \in Q} V_{x}\left( {x_\ast } \right) \leq \bar{R}$.

Stochastic Mirror Descent \cite{2,3} iterates as follows, starting from $x^0 \in Q$,
\begin{equation}
\label{eq2}
x^{k+1}={\rm Mirr}_{x^k} \left( {h\nabla _x f\left( {x^k,\xi ^k} 
\right)} \right), \quad
{\rm Mirr}_{x^k} \left( v \right):=\arg \mathop {\min }\limits_{x\in 
Q} \left\{ {\left\langle {v,x-x^k} \right\rangle +V_{x^k}\left( {x} 
\right)} \right\},
\end{equation}
where $h >0 $ is the stepsize,  $\left\{ {\xi ^k} \right\}_{k \geq 0}$ is an i.i.d. sample of $\xi$.
The main property of SMD-step \cite{3} is
\begin{equation}
\notag
2V_{x^{k+1}}\left( {x} \right)\le \;2V_{x^k}\left( {x} \right)+2h\left\langle 
{\nabla _x f\left( {x^k,\xi ^k} \right),x-x^k} \right\rangle +h^2\left\| 
{\nabla _x f\left( {x^k,\xi ^k} \right)} \right\|_\ast ^2, \quad \forall x \in Q .
\end{equation}
Further, using convexity of $f(x)$, for any $\nabla f(x^k) \in \partial f(x^k)$ and $x \in Q$,  
\[
f\left( {x^k} \right)-f\left( x \right)\le \left\langle {\nabla f(x^k),x^k-x} \right\rangle \le \left\langle {\nabla f(x^k)-\nabla _x f\left( {x^k,\xi ^k} \right),x^k-x} \right\rangle +
\]
\begin{equation}
\notag
+\frac{1}{h}\;\left( {V_{x^k}\left( {x} \right)-V_{x^{k+1}}\left( {x} \right)} 
\right)+\frac{h}{2}\left\| {\nabla _x f\left( {x^k,\xi ^k} \right)} 
\right\|_\ast ^2.
\end{equation}
Taking conditional expectation w.r.t. $\xi ^1,...,\xi ^{k-1}$ and using \eqref{eq:As}, one obtains	  
\[
f\left( {x^k} \right)-f\left( x \right)\le 
 \frac{1}{h}\left( {V_{x^k}\left( {x} \right)-\E\left[ {\left. V_{x^{k+1}}\left( {x} \right) \right|\xi ^1,...,\xi ^{k-1}} \right]} 
\right)+\frac{h}{2}\underbrace {\E\left[ {\left. {\left\| {\nabla _x 
f\left( {x^k,\xi ^k} \right)} \right\|_\ast ^2 } \right|\xi ^1,...,\xi 
^{k-1}} \right]}_{\stackrel{\eqref{eq:As}}{\leq} M^2}.
\]
Since $\left\{ {\xi ^k} \right\}_{k \geq 0}$ is an i.i.d. sample, taking full expectation from the both sides of these inequalities for $k=0,...,N-1$, summing them up and taking $x=x_\ast $ gives, by convexity of $f\left( x \right)$,
\begin{equation}
\notag
\E\left[ {f\left( {\bar {x}^N} \right)} \right]-f_\ast \le \E \frac{1}{N}\sum\limits_{k=0}^{N-1} {f(x^k)} - f_*\leq   \frac{1}{hN}V_{x^0}\left( {x_\ast } \right)+\frac{M^2h}{2} \le \sqrt 
{\frac{2M^2R^2}{N}} ,
\end{equation}
where 
\begin{equation}
R \quad \text{is s.t.} \quad V_{x^0}\left( {x_\ast } \right) \leq R^2,
\quad
\bar {x}^N:=\frac{1}{N}\sum\limits_{k=0}^{N-1} {x^k} ,
\quad
h=\frac{R}{M}\sqrt {\frac{2}{N}}. 
\label{eq:MD_param}
\end{equation}
Choosing 
\begin{equation}
\label{eq6}
N=\left\lceil \frac{2M^2R^2}{\varepsilon ^2}\right\rceil,
\end{equation}
we obtain that $\bar {x}^N$ satisfies $\E\left[ {f\left( {\bar {x}^N} \right)} \right]-f_\ast \le \e$ and, hence, is an $\e$-solution.
Note that this bound for $N$ is optimal \cite{3} up to a multiplicative constant factor for the 
class of convex stochastic programming problems \eqref{eq1} with a.e. bounded stochastic subgradients. 

It turns out that it is possible to prove bounds for probability of large deviation for $f\left( {\bar {x}^N} \right)-f_\ast$.
\begin{Prop}[\cite{2,5,6}]
\label{Prop:1}
Assume that one of the following assumptions holds. a) $\left\| {\nabla _x f\left( {x,\xi } \right)} \right\|_\ast 
\le M$ for a.e. $\xi $; b) $\E_{\xi}\left( {\exp \left( {{\left\| {\nabla _x f\left( {x,\xi } \right)} 
\right\|_\ast ^2 } \mathord{\left/ {\vphantom {{\left\| {\nabla _x f\left( 
{x,\xi } \right)} \right\|_\ast ^2 } {M^2}}} \right. 
\kern-\nulldelimiterspace} {M^2}} \right)} \right)\le \exp \left( 1 \right)$ and $\ln \sigma ^{-1}\ll N$; c) There exists some $\alpha >2$ s.t., for all, $t \geq 0$, $
\Prob\left( {\frac{\left\| {\nabla f\left( {x,\xi } \right)} \right\|_2^2 
}{M^2}\ge t} \right)\le \frac{1}{\left( {t+1} \right)^\alpha } $,
and $\sigma ^{-1 / \left( {\alpha -1} \right)}\ll N$. Then the point $\bar{x}^N$ generated by SMD \eqref{eq2}, \eqref{eq:MD_param} satisfies
\begin{equation}
\label{eq9}
\Prob\left\{ f\left( {\bar {x}^N} \right)-f_\ast \le \frac{C_1M}{\sqrt N 
}\left( R+C_2\bar{R}\sqrt{\ln \left( 1 / \sigma  \right) } \right) 
\right\}\ge 1-\sigma,
\end{equation}
where in the case a) $C_1 = \sqrt{2}$, $C_2 = 2\sqrt{2}$; in the case b) $C_1 = C_2 = 2\sqrt{2}$; in the case c) $C_1 = C_1 (\alpha)$, $C_2 = 1$.
\end{Prop}

\begin{Cor}
\label{Cor:1}
Let any of three assumptions of Proposition \ref{Prop:1} hold. Choose $N = \left\lceil \frac{CM^2\bar {R}^2}{\varepsilon ^2}\right\rceil$, where the constant $C$ depends on $C_1, C_2$. Then the point $\bar{x}^N$ generated by SMD \eqref{eq2}, \eqref{eq:MD_param} satisfies, for any $c \geq 0$,
\begin{equation}
\label{eq12}
\Prob\left\{ {f\left( {\bar {x}^N} \right)-f_\ast \ge c } \right)\le \Prob\left\{ {\eta \ge 
c} \right\},
\end{equation}
where $\eta \in {\rm N}\left( {\varepsilon ,\varepsilon ^2} \right)$ -- normal random variable with mean $\e$ and variance $\e^2$.
\end{Cor}

\section{Parallelization and bounds for probability of large deviation}
In this section, we first discuss a known way to obtain an $(\e,\sigma)$-solution using a number of $\e$-solutions calculated in parallel. Then, we suggest a new way of doing this without calculation of the objective $f(x)$ value, state and prove the main result.

Assume that $\bar{x}^N$ is an $\e/2$-solution, obtained by SMD \eqref{eq2}, \eqref{eq:MD_param} with $N=\left\lceil \frac{8M^2R^2}{\varepsilon ^2}\right\rceil$.
Then, using the Markov inequality \cite{5}, we obtain 
\[
\Prob\left( {f\left( {\bar {x}^N} \right)-f_\ast \ge \varepsilon } \right)\le 
\frac{\E\left[ {f\left( {\bar {x}^N} \right)} \right]-f_\ast }{\varepsilon 
}\le \frac{1}{2}.
\]
If one calculates in parallel $K=\left\lceil \log _2 \left( {\sigma ^{-1}} \right)\right\rceil$ 
independent SMD $\e/2$-solutions $\left\{ {\bar {x}^{N,i}} \right\}_{i=1}^K$ and chooses the one $\bar {x}_{\min }^N $ which minimizes $f\left( {\bar {x}^{N,i}} \right)$, then, in total $\left\lceil \frac{8M^2R^2}{\varepsilon ^2}\right\rceil \left\lceil \log _2 \left( {\sigma ^{-1}} \right)\right\rceil$ calculations of stochastic subgradient, one obtains $\Prob\left( {f\left( {\bar {x}_{\min }^N } \right)-f_\ast \ge \varepsilon } \right)\le \sigma$. Thus, $\bar {x}_{\min }^N $ is an $(\e,\sigma)$-solution of \eqref{eq1}. The crucial point here is the possibility to calculate the value of the function $f\left( x \right)$.

We now suggest a technique which does not rely on the assumption of the function $f\left( x \right)$ value availability. This assumption may not hold \cite{1} in many real stochastic programming problems, e.g. in maximum likelihood approach used in
mathematical statistics.


\begin{Th}
Let any of three assumptions of Proposition \ref{Prop:1} hold. Let $K=\left\lceil 2 \ln \left( {\sigma ^{-1}} \right)\right\rceil$ and $\left\{ {\bar {x}^{N,i}} \right\}_{i=1}^K$ be independent points obtained by SMD \eqref{eq2}, \eqref{eq:MD_param} with $N = \left\lceil \frac{4CM^2\bar {R}^2}{\varepsilon ^2}\right\rceil$. Then the point $\bar 
{x}^K=\frac{1}{K}\sum\limits_{i=1}^K {\bar {x}^{N,i}} $ is an $(\e,\sigma)$-solution of \eqref{eq1}.
\end{Th}
\myproof Using Corollary \ref{Cor:1} with $c=\e$, we obtain, for all $i=1,...,K$, 
\begin{equation}
\Prob\left\{ {f\left( {\bar {x}^{N,i}} \right)-f_\ast \ge \e} \right\}\le \Prob\left\{ 
{\eta _i \ge \e } \right\},
\quad
\eta _i \in {\rm N}\left( {\frac{\e}{2},\frac{\e^2}{4}} \right).
\label{eq:Th1_pr}
\end{equation}
By convexity of $f(x)$, since $\left\{ {\bar {x}^{N,i}} \right\}_{i=1}^K $ are i.i.d., $\left\{ {\eta _i } \right\}_{i=1}^K $ are i.i.d.,  
\begin{align}
& \Prob\left\{ {f\left( {\bar {x}^K} \right)-f_\ast \ge \e } \right\}\le 
\Prob\left\{ {\frac{1}{K}\sum\limits_{k=i}^K {\left( {f\left( {\bar {x}^{N,i}} 
\right)-f_\ast } \right)} \ge \varepsilon } \right\}\stackrel{\eqref{eq:Th1_pr}}{\leq}  \\
&\Prob\left\{ {\frac{1}{K}\sum\limits_{i=1}^K {\eta _i } \ge \frac{\e}{2} 
+ \frac{\e}{2}} \right) = \Prob\left\{ \left(\frac{1}{K}\sum\limits_{i=1}^K { 
{\eta _i  } } \right) -\frac{\e}{2} \ge \frac{\e}{2}  \right\}=\Prob\left\{ {\eta 
-\frac{\e}{2} \ge \frac{\e}{2} \sqrt K } \right\} \leq \sigma ,
\end{align}
where $\eta \in {\rm N}\left( {\frac{\e}{2} ,\frac{\e^2}{4}} \right)$. Here we used well-known facts about the 
properties of sum of independent normal random variables  \cite{4}
\[
{\rm N}\left( \frac{\e}{2} ,\frac{\e^2}{4} \right)+...+{\rm N}\left( 
\frac{\e}{2} ,\frac{\e^2}{4} \right)\mathop =\limits^d {\rm N}\left( 
\frac{K\e}{2} ,\frac{K\e^2}{4} \right),
\quad
\frac{1}{K}{\rm N}\left( \frac{K\e}{2} ,\frac{K\e^2}{4} \right) \mathop =\limits^d {\rm 
N}\left( {\frac{\e}{2},\frac{\e^2}{4K}} \right),
\] 
that $K=\left\lceil2 \ln \left( {\sigma ^{-1}} \right)\right\rceil$ and, for $\eta \in {\rm N}\left( {\frac{\e}{2} ,\frac{\e^2}{4}} \right)$, $\Prob\left\{ {\eta  -\frac{\e}{2} \ge \frac{\e}{2} \sqrt{2\ln \left( {\sigma ^{-1}} \right)} } \right\} \leq \sigma$.
Thus, the point $\bar {x}^K$ is an $(\e,\sigma)$-solution of \eqref{eq1}. \qed

To sum up, under some mild assumptions, but without calculation of the objective function value, we propose a way to obtain an $(\e,\sigma)$-solution using a number of $\e$-solutions calculated in parallel. This approach allows to reduce the computational time without any loss in the solution quality. At the same time, we answer to the question of Yu. Nesterov \cite{7}. The question can be stated as follows. When the quality of a solution of problem \eqref{eq1} obtained by one wise old man, thinking for $\Theta\left({M^2\bar {R}^2\ln \left( {\sigma 
^{-1}} \right)} /\e^2 \right)$ days, is the same as obtained by $\Theta \left( \ln \left( {\sigma 
^{-1}} \right)\right)$ experts, each thinking for $\Theta\left( M^2\bar {R}^2/ \e^2 \right)$ days? Our answer is that the quality is equivalent under any of three assumptions of Proposition \ref{Prop:1}.




\begin{thebibliography}{9}

\bibitem{1} \textit{Shapiro A., Dentcheva D., Ruszczynski A.} Lecture on stochastic programming. Modeling and theory. -- MPS-SIAM series on Optimization, 2014.
\bibitem{2} \textit{Nemirovski A.,  Juditsky A., Lan, G., Shapiro, A.   } Robust  stochastic  approximation approach to stochastic programming. SIAM J. Optim. 2009, V. 19. P.1574-1609.
\bibitem{3} \textit{Nemirovski A.} Lectures on modern convex optimization analysis, algorithms, and engineering applications. Philadelphia: SIAM, 2013.
\bibitem{4} \textit{Durrett R.} Probability and examples. -- Cambridge University Press, 2010.
\bibitem{5} \textit{Guiges V., Juditsky A., Nemirovski A. }Non-asymptotic confidence bounds for the optimal value of a stochastic program // e-print, 2016 arXiv:1601.07592
\bibitem{6} \textit{Gasnikov A.V. }Searching equilibriums in large transport networks. Doctoral thesis. -- MIPT, 2016. [in Russian] arXiv:1607.03142
\bibitem{7} \textit{Nesterov Yu., Vial J.-Ph.} Confidence level solution for stochastic programming // Automatica. 2008. V. 44. no. 6. P. 1559--1568.
\end{thebibliography}
\end{document}